\documentclass[11pt]{article}
\pdfoutput=1                            % to force arxiv to use pdflatex

\usepackage[utf8]{inputenc}				% handle accents properly
\usepackage[english]{babel}				% handle hyphenation properly
\usepackage[hidelinks]{hyperref}		% hyperlinks

\usepackage{fullpage}					% use the whole page
\usepackage[shortlabels]{enumitem}		% customizable itemization
\usepackage[compress]{cite}				% citation shortening

\usepackage{caption}
\usepackage{amsmath,amsfonts,amsthm}
\usepackage{booktabs}
\newtheorem{thm}{Theorem}

\usepackage{tikz}
\usetikzlibrary{calc}
\usepackage{pgfplots}
\pgfplotsset{compat=1.17}
\usepgfplotslibrary{fillbetween}
\graphicspath{{./figures/}}

\newcommand{\norm}[1]{\|#1\|}
\newcommand{\grad}{\nabla}
\newcommand{\R}{\mathbb{R}}
\newcommand{\tp}{\mathsf{T}}
\newcommand{\bmat}[1]{\begin{bmatrix}#1\end{bmatrix}}
\DeclareMathOperator*{\argmin}{\arg\min}

\DeclareMathOperator{\prox}{prox}
\newcommand{\suchthat}{\;\ifnum\currentgrouptype=16 \middle\fi|\;}
\newcommand{\set}[2]{\left\{#1\suchthat#2\right\}}
\usepackage{colonequals}
\newcommand{\defeq}{\colonequals}

\let\OLDthebibliography\thebibliography
\renewcommand\thebibliography[1]{
  \OLDthebibliography{#1}
  \setlength{\parskip}{0pt}
  \setlength{\itemsep}{4pt plus 0.3ex}
}

\title{The Analysis of Optimization Algorithms\\
\Large A dissipativity approach}
\author{Laurent Lessard\\
Department of Mechanical and Industrial Engineering\\
Northeastern University ({\small\texttt{l.lessard@northeastern.edu}})}
\date{}

\begin{document}
\maketitle

Optimization algorithms are ubiquitous across all branches of engineering, computer science, and applied mathematics. Typically, such algorithms are iterative in nature and seek to approximate the solution to a difficult optimization problem. The general form of an optimization problem is to
\begin{align*}
	\text{minimize}&\qquad f(x) \\
	\text{subject to}&\qquad x \in C,
\end{align*}
where $x \in \R^d$ is the set of $d$ real \emph{decision variables}, $C$ is the \emph{feasible set}, and $f:\R^d\to \R$ is the \emph{objective function}. The goal is to find $x\in C$ such that $f(x)$ is as small as possible.
For example, in structural mechanics, $x$ could be lengths and widths of the members in a truss design, $C$ could encode stress limitations of the materials and load bearing constraints, and $f$ could be the total cost of the design. Solving this optimization problem would provide the cheapest truss design that satisfies all design specifications.
In statistics, $x$ could be parameters in a statistical model, $C$ could correspond to constraints such as certain parameters being positive, and $f$ could be the negative log-likelihood of the model given the observed data. Solving this optimization problem finds the maximum likelihood model.
While some optimization problems can be solved analytically (such as least squares problems), analytical solutions do not exist for most optimization models, and numerical methods must be used to approximate the solution. Even when an analytical solution exists, it may be preferable to use numerical methods because they can be more computationally efficient.
Iterative optimization algorithms typically begin with an estimate $x^0$ of the solution, and each step refines the estimate, producing a sequence $x^0, x^1, \dots$. If properly designed, then $x^k \to x^\star$ in the limit, and as many iterations as needed can be used to achieve the desired level of accuracy.
Depending on the nature and structure of $f$ and $C$ in the optimization problem, different optimization algorithms may be appropriate. The design, selection, and tuning of optimization algorithms is more of an art than a science. Many different algorithms have been proposed, some with theoretical guarantees and others with a strong record of empirical performance. In practice, some algorithms converge slowly yet are predictable and reliable. Meanwhile, other algorithms converge more rapidly on average, but can fail spectacularly in some cases. Algorithm selection and tuning is typically performed by experts with deep area knowledge.
In many ways, iterative algorithms behave like control systems, and the choice of algorithm is akin to the choice of controller. In the sections that follow, we formalize this connection and describe how optimization algorithms can be viewed as controllers performing robust control. This work will also show how dissipativity theory can be used to analyze the performance of many classes of optimization algorithms. This allows selection and tuning of optimization algorithms to be performed in an automated and systematic way.
%%%%%%%%%%%%%%%%%%%%%%%%%%%%%%%%%%%%%%%%%%%%%%%%
\section{Black-box paradigm and performance evaluation}
To reason about different optimization algorithms and their performance, it is common to employ a \emph{black-box} paradigm~\cite[\S1.1.2]{nesterov}. This model assumes the availability of \emph{oracles} that can be queried to provide pertinent information about the objective function or the constraints. The oracles are how the algorithm interfaces with the optimization problem.
For example, if $f$ is differentiable and $C$ is a convex set, a popular algorithm is \emph{projected gradient descent}, which begins with a some guess value $x^0$ and follows the iterations:
\begin{equation}\label{eq:pgd}
	x^{k+1} = \Pi_C \left( x^k - \eta \grad f( x^k ) \right),
	\qquad\text{for }k=0,1,\dots
\end{equation}
Here, $\eta > 0$ is a tuning parameter called the \emph{stepsize}, 
$\grad f$ is the gradient of $f$, and $\Pi_C$ is the projection onto the set $C$. For sufficiently small $\eta$ and under appropriate regularity conditions on $f$, projected gradient descent will converge to a solution $x^\star$ of the constrained optimization problem. This example includes two oracles:
\begin{enumerate}
	\item Gradient oracle: Given $x$, return $\grad f (x)$.
	\item Projection oracle: Given $x$, return $\Pi_C(x)$.
\end{enumerate}
In the black-box model, each oracle query incurs a fixed cost (usually \emph{time}) and all other costs associated with computer storage or memory are ignored. The \emph{performance} of an iterative method is based on the total time required for an error measure to reach a specified value. Common error measures include distance to optimality $\norm{x^k-x^\star}$ and function error $f(x^k) - f(x^\star)$.
Iterative algorithms typically call each oracle once per iteration, so performance can be evaluated by measuring the \emph{convergence rate}, which is how quickly the error measure decreases per iteration. For example, an algorithm exhibits \emph{geometric} convergence if there is some $\rho \in [0,1)$ and $K > 0$ such that
\begin{equation}\label{def:gemoconv}
\norm{x^k - x^\star} \leq K\rho^k \norm{x^0 - x^\star}\quad\text{for }k=0,1,\dots.
\end{equation}
Smaller $\rho$ corresponds to a faster algorithm. The smallest value of $\rho$ that satisfies~\eqref{def:gemoconv} can depend on the oracles used (the choice of $f$ and $C$) and the initial condition $x^0$. Care must be taken when interpreting convergence rates, because as $\rho$ becomes smaller, $K$ may get larger, and as $\rho$ approaches its minimum value, then we may have $K\to\infty$.
\paragraph{Algorithm analysis}
Algorithm analysis is the practice of determining bounds on the rate of convergence of algorithms subject to various assumptions. Conventionally, algorithm analysis provides an assurance that a given algorithm will work well for many instances of an optimization problem. For example, it may be desirable for algorithm $\mathcal{A}$ to converge rapidly for many different pairs of oracles $(f,C)$. We call this set of admissible oracle pairs $\mathcal{F}$.
A typical analysis query might be: What is the \emph{worst-case} geometric convergence rate $\rho$ achieved by $\mathcal{A}$ over the set $\mathcal{F}$? Mathematically, this is given by the expression
\begin{equation}\label{def:conv_rate}
	\rho(\mathcal{A},\mathcal{F}) \defeq \inf \set{ \rho > 0 }{\sup_{(f,C)\in\mathcal{F}}\,\sup_{x^0}\, \sup_{k\ge 0}\, \frac{\norm{x^k-x^\star}}{\rho^k \norm{x^0-x^\star}}<\infty }.
\end{equation}
Equation~\eqref{def:conv_rate} states that for any $\rho > \rho(\mathcal{A},\mathcal{F})$, the geometric convergence criterion~\eqref{def:gemoconv} holds for all choices of oracles $(f,C)\in\mathcal{F}$ and all initial conditions $x^0$. So, $\rho(\mathcal{A},\mathcal{F})$ is the fastest convergence rate that is guaranteed to hold over all admissible problem instances and initial conditions. Two algorithms $\mathcal{A}_1$ and $\mathcal{A}_2$ can then be compared based on their convergence rates. If $\rho(\mathcal{A}_1,\mathcal{F}) < \rho(\mathcal{A}_2,\mathcal{F})$, then $\mathcal{A}_1$ is faster than $\mathcal{A}_2$ in the worst case.
The notion of \emph{worst-case convergence rate} in algorithm analysis is akin to the notion of \emph{robust stability} in nonlinear control. The following sections explore this connection in greater detail and show how tools from robust control can be brought to bear on the problem of algorithm analysis.
%%%%%%%%%%%%%%%%%%%%%%%%%%%%%%%%%%%%%%%%%%%%%%%%
\section{Algorithm analysis as robust control}
The performance evaluation of iterative algorithms under the black-box paradigm may be reframed as certifying robust stability for a feedback system. Specifically, the algorithm can be written as a discrete-time dynamical system in feedback with its oracles. It will be illustrated through examples how a variety of algorithms can be converted to \emph{feedback form}.
\subsection{Projected gradient descent}
\begin{figure}[ht]
	\centering
	\includegraphics{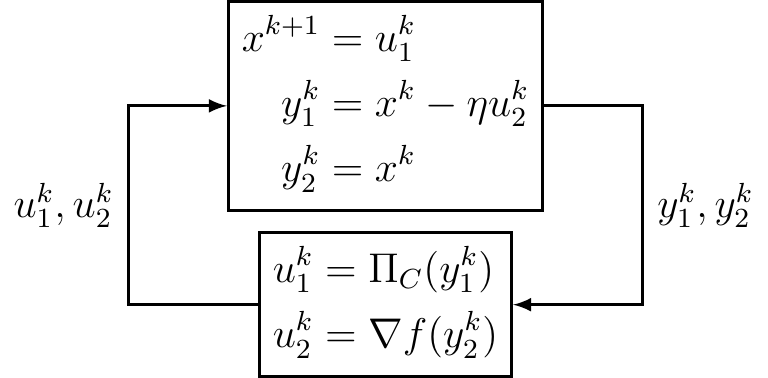}
	\caption{Feedback interconnection for projected gradient descent, an iterative algorithm with update equations given by~\eqref{eq:pgd}. The feedback form separates the algorithm dynamics (which are a linear time-invariant (LTI) system) from the oracle calls (which are treated as unknown nonlinearities).}
	\label{fig:blkdiag_projgrad}
\end{figure}

Returning to the projected gradient descent example~\eqref{eq:pgd}, the algorithm has access to two oracles: a projection operator $\Pi_C$ and the gradient $\grad f$. If the problem is unconstrained, the simplification $\Pi_C = I$ occurs and~\eqref{eq:pgd} becomes ordinary gradient descent.
It is converted to feedback form by defining the input-output pairs of $\Pi_C$ and $\grad f$ as $(y_1,u_1)$ and $(y_2,u_2)$, respectively. The ensuing block diagram is illustrated in Figure~\ref{fig:blkdiag_projgrad}.
\subsection{Nesterov's accelerated method}
\begin{figure}[ht]
	\centering
	\includegraphics{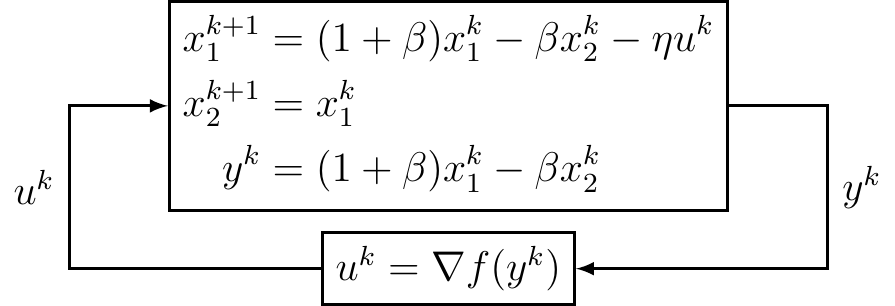}
	\caption{Feedback interconnection for Nesterov's accelerated method, whose update equations are given in~\eqref{eq:nesterovupdate}.}
	\label{fig:blkdiag_nesterov}
\end{figure}

Nesterov's accelerated method~\cite[\S2.2]{nesterov}  is a popular iterative approach used to solve the unconstrained optimization problem $\min_x\,f(x)$, where $f$ is continuously differentiable, and access to a gradient oracle $\grad f$ is provided. The algorithm has two states $(x^k,y^k)$ and uses the update
\begin{subequations}\label{eq:nesterovupdate}
\begin{align}
	y^{k} &= x^{k} + \beta (x^{k} - x^{k-1}) \label{eq:nestup1}\\
	x^{k+1} &= y^k - \eta \grad f(y^k) \label{eq:nestup2}.
\end{align}
\end{subequations}
The state $y^k$ extrapolates based on the current iterate $x^k$ and previous iterate $x^{k-1}$. Then, gradient descent is performed based on $y^k$. The \emph{momentum parameter} $\beta$ controls the amount of extrapolation. When $\beta = 0$, $y^k=x^k$ and the familiar gradient descent algorithm is recovered. The idea is to tune $\beta$ to obtain faster convergence than ordinary gradient descent.
To convert~\eqref{eq:nesterovupdate} to feedback form, substitute~\eqref{eq:nestup1} into~\eqref{eq:nestup2} and define the new state variables $x_1^k = x^k$ and $x_2^k = x^{k-1}$. We then obtain the feedback form illustrated in Figure~\ref{fig:blkdiag_nesterov}.
\subsection{ADMM algorithm}
\begin{figure}[ht]
	\centering
	\includegraphics{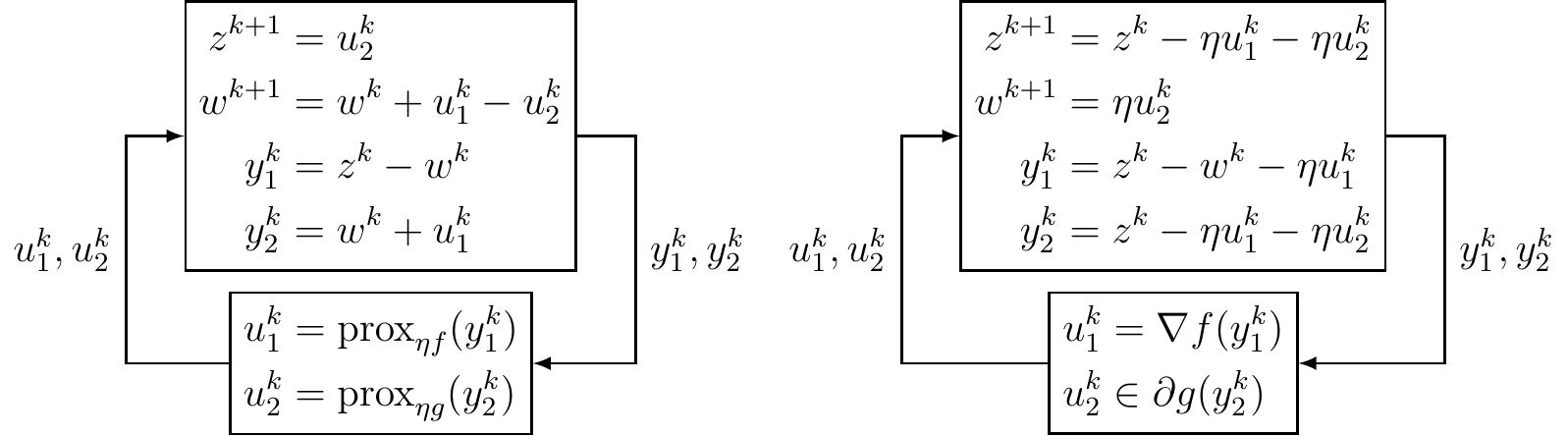}
	\caption{Equivalent feedback interconnections for the Alternating Direction Method of Multipliers (ADMM) applied to composite unconstrained optimization. \textbf{Left:} An explicit loop that uses proximal operators as oracles~\eqref{eq:admm_prox}. \textbf{Right:} An implicit loop that uses gradient and subgradient oracles~\eqref{eq:admm_grad}. Other representations are possible, for example using $\grad f$ and $\prox_{\eta g}$. Any of these representations of ADMM can be used for analysis in the dissipativity framework and yield the same results.}
	\label{fig:blkdiag_admm}
\end{figure}

The Alternating Direction Method of Multipliers (ADMM)~\cite{admm} is a popular algorithm used to solve \emph{composite} optimization problems. That is, the objective function can be split into two parts with different properties. The canonical problem takes the form:
\begin{align*}
	\text{minimize} \qquad & f(x) + g(z) \\
	\text{subject to:}\qquad& Ax + Bz = c.
\end{align*}
For example, $f$ may be convex and differentiable, while $g$ may be a nondifferentiable regularization term or the indicator set for a convex constraint. The ADMM algorithm has three state variables $(x^k,z^k,w^k)$ and a tuning parameter $\eta$. The algorithm uses the update
\begin{align*}
x^{k+1} &= \argmin_x\, f(x) + \tfrac{1}{2\eta}\norm{ Ax + Bz^k - c + w^k}^2, \\
z^{k+1} &= \argmin_z\, g(z) + \tfrac{1}{2\eta}\norm{ Ax^{k+1} + Bz - c+ w^k}^2, \\
w^{k+1} &= w^k + A x^{k+1} + B z^{k+1} - c.
\end{align*}
To simplify exposition, consider the composite unconstrained case: $\min_x f(x)+g(x)$. This is the special case with $A=I$, $B=-I$, and $c=0$. The ADMM algorithm can be written in feedback form in several equivalent ways, depending on which oracles are used. For example, the \emph{proximal operator}~\cite{prox} is defined as $\prox_{f}(z) \defeq \argmin_x f(x) + \frac{1}{2}\norm{x-z}^2$. Use it to rewrite the ADMM update equations as
\begin{subequations}\label{eq:admm_prox}
	\begin{align}
		x^{k+1} &= \prox_{\eta f} (z^k - w^k), \\
		z^{k+1} &= \prox_{\eta g} (x^{k+1} + w^k), \\
		w^{k+1} &= w^k + x^{k+1} - z^{k+1}.
	\end{align}
\end{subequations}
Alternatively, if $f$ is differentiable and $g$ is convex but not differentiable, replace the $\prox$ updates by their corresponding first-order optimality conditions. This yields
\begin{subequations}\label{eq:admm_grad}
	\begin{align}
		0 &= \grad f(x^{k+1}) + \tfrac{1}{\eta} (x^{k+1} - z^k + w^k), \\
		0 &\in \partial g(z^{k+1}) + \tfrac{1}{\eta} (z^{k+1} - x^{k+1} - w^k), \\
		w^{k+1} &= w^k + x^{k+1} - z^{k+1},
	\end{align}
\end{subequations}
where $\partial g(z) \defeq \left\{ v \in \R^d \; \mid\; g(x)-g(z) \geq v^\tp (x-z)\text{ for all }z\in \R^d\right\}$ denotes the set of subgradients of $g$. If $g$ is differentiable, then $\partial g(z) = \{\grad g(z)\}$. Both~\eqref{eq:admm_prox} and \eqref{eq:admm_grad} can be put in feedback form to obtain the block diagrams illustrated in Figure~\ref{fig:blkdiag_admm}. In both cases, the state variable $x^k$ can be eliminated, so only two states are needed.
The ADMM representation that uses gradients and subgradients in Figure~\ref{fig:blkdiag_admm} is \emph{implicit} because it contains a circular dependency: $y_1^k$ depends on $u_1^k$, which in turn depends on $y_1^k$. Therefore, this feedback representation cannot be used as a substitute for an implementation such as~\eqref{eq:admm_prox}. Nevertheless, the implicit representation can still be used in dissipativity theory for algorithm analysis.
\subsection{More general algorithms}
In all the cases above, the algorithms can be expressed in the form of a linear time-invariant (LTI) system $G$ in feedback with the oracles. Letting $\xi^k$, $u^k$, and $y^k$ be the concatenated states, oracle outputs, and oracle inputs, the algorithms are represented in the following general form.
\noindent\begin{minipage}[T]{0.45\linewidth}
\begin{subequations}\label{eq:genalg}
\begin{align}
\xi^{k+1} &= A \xi^k + B u^k \\
y^k &= C\xi^k + D u^k \\
u^k = \bmat{u_1^k \\ \vdots \\ u_m^k} &= \bmat{ \phi_1(y_1^k) \\ \vdots \\ \phi_m(y_m^k)} = \phi(y^k).
\end{align}
\end{subequations}
\end{minipage}
\hfill
\begin{minipage}[T]{0.55\linewidth}
	\centering
	\begin{tikzpicture}[thick,>=latex,auto]
		\node[draw] (alg) {$G:\, \left\{\begin{aligned}
			\xi^{k+1} &= A \xi^k + B u^k \\
			y^k &= C \xi^k + Du^k \\
			\end{aligned}\right\}$};
		\node[draw,anchor=north] (oracle) at ($(alg.south)+(0,-0.2)$) {$u^k = \phi(y^k)$};
		\draw[->] (alg.east) -- +(0.7,0) |- node[pos=0.25]{$y^k$} (oracle);
		\draw[<-] (alg.west) -- +(-0.7,0) |- node[pos=0.25,swap]{$u^k$}(oracle);
		\end{tikzpicture}
\end{minipage}

\noindent The fact that $G$ is LTI will be important for our dissipativity analysis, as it will allow us to search for Lyapunov functions in a tractable manner. Not all algorithms have a feedback form with an LTI $G$. For example,
\begin{itemize}
	\item Algorithms with parameters that change on a fixed schedule, such as gradient descent with a diminishing stepsize, will yield a feedback representation where $G$ is a linear time-varying (LTV) system.
	\item Algorithms with parameters that change adaptively, such as the nonlinear conjugate gradient method, will yield a feedback representation where $G$ is a linear parameter-varying (LPV) system.
	\item Algorithms where the state updates are not linear functions of the previous state or oracle outputs will yield a feedback representation where $G$ is nonlinear.
\end{itemize}
Despite these limitations, algorithms can still generally be written as a feedback interconnection of some system $G$ and the set of oracles (nonlinearities). Since dissipativity theory can be applied to any system~\cite{willems1} (including LTV, LPV, and systems with nonlinear dynamics), the dissipativity approach can in principle be used to analyze any iterative algorithm.
%%%%%%%%%%%%%%%%%%%%%%%%%%%%%%%%%%%%%%%%%%%%%%%%
\section{Dissipativity theory}
Dissipativity theory may be viewed as a counterpart to Lyapunov theory but for systems with inputs. Consider a discrete-time dynamical system satisfying the state-space equation
\[
\xi^{k+1} = A \xi^k + B u^k.
\]
In classical dissipativity theory~\cite{willems1,willems2}, $u^k$ is an external supply that drives the dynamics governed by the state $\xi^k$. For example, in a mechanical system, $u^k$ would be a vector of external forces and torques, while $\xi^k$ would be a vector of generalized coordinates such as positions and velocities. The two key concepts in dissipativity are \emph{storage} and \emph{supply}.
\begin{enumerate}
	\item The \emph{storage function} $V(\xi^k)$ can be interpreted as a notion of stored energy. In our mechanical example, this would be the total energy (kinetic and potential) in the system. The storage function always satisfies $V(\xi^k) \geq 0$.
	\item The \emph{supply rate} $S(\xi^k,u^k)$ can be interpreted as a notion of work done by the external forces and torques. The supply rate may also depend on the current values of the generalized coordinates. When $S > 0$, the external force is \emph{adding} energy to the system. When $S < 0$, the external force is \emph{extracting} energy from the system.
\end{enumerate}
A \emph{dissipation inequality} states that the change in stored energy can be no greater than the energy provided by the external supply:
\begin{equation}\label{eq:diss_will}
V(\xi^{k+1}) - V(\xi^k) \leq S(\xi^k, u^k).
\end{equation}
If a dissipation-like inequality holds and the supply rate is of a particular form, useful stability properties of the system can be deduced. We now consider different types of supply rates that are relevant for optimization algorithms.
%%%%%%%%%%%%%%%%%%%%%%%%%%%%%%%%%%%%%%%%%%%%%%%%
\section{Supply rates for families of oracles}
In the context of robust control, the most well-studied classes of nonlinearities include sector-bounded and slope-restricted nonlinearities, because they can be used to model common nonlinear phenomena such as saturation and stiction. In the context of optimization algorithms, oracles can often have similar properties, thus creating parallels between the two application areas. For each type of nonlinearity or oracle, we describe how to formulate an appropriate supply rate that can be used to analyze optimization algorithms.
\subsection{Sector-bounded oracle}
\begin{figure}[ht]
	\centering
	\includegraphics{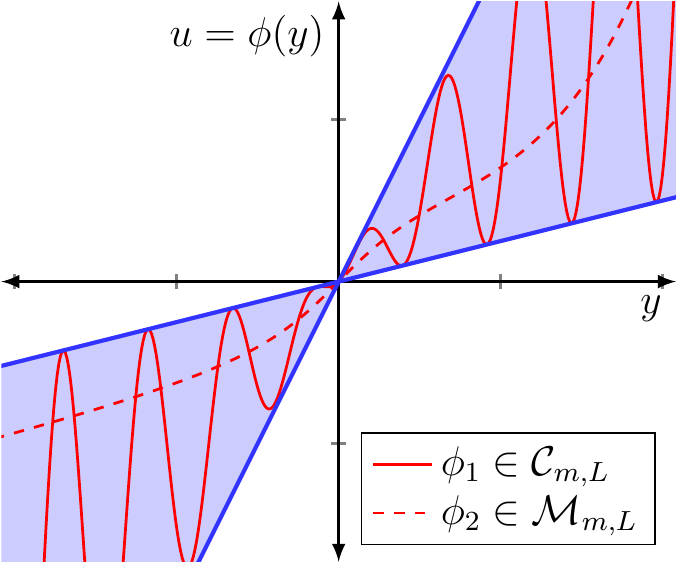}
	\caption{One-dimensional examples of a sector-bounded nonlinearity $\phi_1 \in \mathcal{C}_{m,L}$ and a slope-restricted nonlinearity $\phi_2 \in \mathcal{M}_{m,L}$, both with $m=\tfrac{1}{4}$ and $L = 2$. In the sector-bounded case, all input-output pairs $(u,y)$ satisfy the quadratic inequality~\eqref{eq:sector} (which is $\tfrac{1}{2}y^2 - \tfrac{9}{4}yu + u^2 \leq 0$), shown as the shaded region.
	In the slope-restricted case, a stronger condition is satisfied:
	$\tfrac{1}{2}(y_1-y_2)^2 - \tfrac{9}{4}(y_1-y_2)(u_1-u_2) + (u_1-u_2)^2 \leq 0$, for all $y_1,y_2,u_1,u_2 \in \R$. In other words, the slope of the line connecting any pair of points on the graph must be between $\frac{1}{4}$ and $2$.}
	\label{fig:sector}
\end{figure}

An oracle $u = \phi(y)$ is \emph{sector-bounded} with lower bound $m$ and upper bound $L$ with $m < L$ if the input-output pair $(y,u)$ satisfies 
\begin{equation}\label{eq:sector}
	S_\mathcal{C}(y,u) \defeq 
\bmat{y\\u}^\tp
\bmat{ mL & -\frac{L+m}{2} \\ -\frac{L+m}{2} & 1}
\bmat{y\\u} \leq 0.
\end{equation}
In the general dissipativity framework, the supply rate $S(\xi,u)$ may depend on the state $\xi$ and the input $u$. The form of the supply rate $S_\mathcal{C}$ in~\eqref{eq:sector} still fits into this framework because the oracle is a function of $y$, which will depend on $\xi$ and $u$ through the algorithm update equations~\eqref{eq:genalg}. In this case, $S_\mathcal{C} \leq 0$. Thus, in the energy interpretation of dissipativity, the \emph{external force} is extracting energy from the system.
The inequalities~\eqref{eq:sector} hold \emph{pointwise} in time for any pair $(y,u)$, so the graph of a sector-bounded oracle is contained in the \emph{interior conic} region illustrated in Figure~\ref{fig:sector}. Denote the set of all sector-bounded nonlinearities with parameters $(m,L)$ with the symbol $\mathcal{C}_{m,L}$. A useful pair of inequalities that follow from~\eqref{eq:sector} are given by
\begin{gather}\label{eq:sector2}
	m \norm{y} \leq \norm{u} \leq L \norm{y} \\
	m \norm{y}^2 \leq u^\tp y \leq L \norm{y}^2.
	\label{eq:sector3}
\end{gather}
See \ref{sb:sproc} for insight on how inequalities such as~\eqref{eq:sector2}--\eqref{eq:sector3} can be proved.
Sector-bounded nonlinearities are widely studied in controls, dating back to the introduction of \emph{absolute stability} by Lur'e and Postnikov~\cite{lure_postnikov}. In this setting, $\phi$ is a static nonlinearity. Special cases include the small-gain theorem and passivity theory~\cite{zames}, where $\norm{\phi(y)} \leq \gamma \norm{y}$ ($m=-\gamma$ and $L = \gamma$) and $y^\tp \phi(y) \geq 0$ ($m=0$ and $L \to\infty$), respectively.
Sector-bounded oracles can occur in optimization when oracles are subject to multiplicative noise. For example, round-off error may occur due to finite-precision arithmetic. Alternatively, a call to an oracle may involve a complicated simulation that inherently produces approximate results due to time budget  limitations.
If $y$ is the true signal, multiplicative noise transforms the signal into $u = y + \delta_y$, where $\norm{\delta_y} \leq \varepsilon \norm{y}$. Here, $\varepsilon$ is the noise strength ($\varepsilon = 0.1$ would correspond to 10\% multiplicative noise). Rearranging this inequality yields $\norm{y-u} \leq \varepsilon \norm{y}$. Comparing to~\eqref{eq:sector}, this corresponds to a sector-bounded nonlinearity with $m=1-\varepsilon$ and $L=1+\varepsilon$.
\subsection{Slope-restricted oracles}
An oracle $u = \phi(y)$ is \emph{slope-restricted} with lower bound $m$ and upper bound $L$ with $m < L$ if every pair of input-output pairs $(y_1,u_1)$ and $(y_2,u_2)$ satisfies
\begin{equation}\label{eq:monotone}
	S_{\mathcal{M}}(y_1,y_2,u_1,u_2) \defeq \bmat{y_2-y_1\\u_2-u_1}^\tp
	\bmat{ mL & -\frac{L+m}{2} \\ -\frac{L+m}{2} & 1}
	\bmat{y_2-y_1\\u_2-u_1} \leq 0.
\end{equation}
In a manner analogous to how~\eqref{eq:sector2} and~\eqref{eq:sector3} were derived for sector-bounded nonlinearities, it can be shown that slope-restricted nonlinearities enjoy the properties
\begin{gather}\label{eq:monotone2}
	m\norm{y_2-y_1} \leq \norm{u_2-u_1} \leq L\norm{y_2-y_1} \\
	\label{eq:monotone3}
	m\norm{y_2-y_1}^2 \leq (u_2-u_1)^\tp (y_2-y_1) \leq L\norm{y_2-y_1}^2.
\end{gather}
Slope-restricted nonlinearities (also called \emph{incremental}) were studied by Zames and Falb~\cite{zames1968stability} (for example, \emph{incremental passivity} or \emph{incremental small-gain} ). Examples of slope-restricted nonlinearities include saturation or elements exhibiting hysteresis. See Figure~\ref{fig:sector} for a visual example.
Denote the set of all slope-restricted nonlinearities with parameters $(m,L)$ with the symbol $\mathcal{M}_{m,L}$.
In optimization, different types of slope-restrictedness bear different names. The proximal operator $\prox_g$ for any convex function $g$, used for example in ADMM~\eqref{eq:admm_prox}, is \emph{firmly nonexpansive}. That is,
$(x-y)^\tp(\prox_g(x)-\prox_g(y)) \geq \norm{\prox_g(x)-\prox_g(y)}^2$ for all $x,y\in\R^d$~\cite[Prop.~4.16 and~12.28]{bauschke2017convex}. In other words, $\prox_f$ satisfies~\eqref{eq:monotone} with $m=0$ and $L=1$. A special case is when $g$ is the indicator function of a convex set $C$ [that is $g(x) = 0$ if $x\in C$ and $g(x) = +\infty$ otherwise], then $\prox_g(x) = \Pi_C(x)$ is the Euclidean projection of $x$ onto the set $C$.
Another example is the subgradient of a convex function $g$, which were used in~\eqref{eq:admm_grad}. Subgradients are \emph{monotone}, that is, $(x_1 - x_2)^\tp (v_1 - v_2) \geq 0$ for all $x_i\in\R^d$ and $v_i \in \partial g(x_i)$. In other words, $\partial g$ satisfies~\eqref{eq:monotone} with $m=0$ and $L\to\infty$. A special case is when $g$ is quadratic and positive semidefinite: $g(x) = x^\tp A x$, where $A+A^\tp \succeq 0$.
 
\subsection{Gradient of a convex function}
The case where the nonlinearity is the gradient of a convex and continuously differentiable function is of particular interest in the field of optimization, even if it is a less common occurrence in controls. If $f$ is a continuously differentiable function, define the following notions.
\begin{enumerate}
	\item $f$ is \emph{strongly convex} with parameter $m > 0$ if $g(x) \defeq f(x) - \tfrac{m}{2}\norm{x}^2$ is a convex function. In other words, $\theta g(x) + (1-\theta)g(y) \geq g(\theta x + (1-\theta)y)$ for all $x,y\in\R^d$ and $\theta \in [0,1]$.
	\item $f$ has \emph{Lipschitz gradients} with parameter $L > 0$ if $\norm{\grad f(x)-\grad f(y)} \leq L \norm{x-y}$ for all $x,y\in\R^d$.
\end{enumerate}
Strong convexity means that $f$ is not too flat, while Lipschitz gradients ensure that $f$ does not grow too quickly. We now state a useful property of functions that are both strongly convex and have Lipschitz gradients. Denote the set of all nonlinearities satisfying the two properties above with the symbol $\mathcal{F}_{m,L}$.
\begin{thm}
Suppose $f:\R^d\to\R$ is continuously differentiable, strongly convex with parameter $m$, and has Lipschitz gradients with parameter $L$.
For all $y_i \in \R^d$ with $f_i = f(y_i)$ and $u_i = \grad f(y_i)$ for $i=1,2$, 
\begin{multline}\label{eq:strcvx}
	\mathcal{S}_{\mathcal{F}}(y_1,y_2,u_1,u_2) \defeq
	\frac{1}{2(L-m)}
	\bmat{y_2-y_1\\u_2-u_1}^\tp
	\bmat{ mL & -\frac{L+m}{2} \\ -\frac{L+m}{2} & 1}
	\bmat{y_2-y_1\\u_2-u_1}\\ + \frac{1}{2}(u_1+u_2)^\tp (y_2-y_1)
	\leq f_2-f_1.
\end{multline}
\end{thm}
This result is proven for example in~\cite[Thm.~4]{taylor2017smooth}.
In large-scale optimization problems, it is typical to assume the availability of a gradient oracle. Algorithms that use such an oracle are called \emph{first-order methods}. Popular examples of optimization problems that are often solved on a large scale include logistic regression and regularized least-squares problems such as the lasso~\cite{james2013introduction}.
\subsection{Other oracle types}
Many optimization problems involve objective functions (oracles) that are convex, but not strongly convex. This has motivated the study of conditions that are weaker than strong convexity but can still provide useful convergence guarantees. One such example is the \emph{Polyak--{\L}ojasiewicz} condition: 
$\frac{1}{2m}\norm{u}^2 \geq f-f^\star$.
Others include the \emph{error bound} condition, the \emph{quadratic growth} condition, and the \emph{restricted secant inequality}. For a survey of such conditions, refer to~\cite{bolte2017error,karimi2016linear} and references therein.
Although we restrict attention to strong convexity in the present article, the aforementioned alternatives to strong convexity can be used just as easily. All of these conditions are characterized by inequalities similar to \eqref{eq:strcvx} that can be used directly as supply rates; they are quadratic in $(x^k-x^\star)$ and $u^k$, and they are linear in $f^k-f^\star$.
\subsection{Nestedness and supply rates}
The three classes of nonlinearities  characterized by~\eqref{eq:sector}, \eqref{eq:monotone}, and \eqref{eq:strcvx} are nested. That is, $\mathcal{F}_{m,L} \subseteq \mathcal{M}_{m,L} \subseteq \mathcal{C}_{m,L}$.% (see Figure~\ref{fig:venn_diagram}).
This follows because if \eqref{eq:strcvx} is added to itself with indices interchanged, then~\eqref{eq:monotone} is recovered. So, \eqref{eq:monotone} is a special case of \eqref{eq:strcvx}. Moreover, if $y_1 = 0$ (and $u_1=0$) in \eqref{eq:monotone}, then ~\eqref{eq:sector} is recovered. So, \eqref{eq:sector} is a special case of \eqref{eq:monotone}. In the one-dimensional case ($d=1$), all gradients of convex functions are slope-restricted, and vice versa. In other words, $\mathcal{F}_{m,L} = \mathcal{M}_{m,L}$. This is not the case when $d\geq 2$, which is related to the notion of cyclic monotonicity~\cite{rockafellar2015convex}.
The nestedness property implies that a supply rate that is valid for one class of oracles is also valid for any oracle belonging to a subclass.
For each class of oracles, the fundamental inequalities \eqref{eq:sector}, \eqref{eq:monotone}, and \eqref{eq:strcvx} can be directly used as supply rates. Typically, there will be many valid choices of supply rates, and different choices could yield different convergence rate guarantees. The case studies that follow show how to optimize the choice of supply rate to produce the least conservative estimate of convergence rate attainable.
% \begin{figure}[ht]
% \centering
% \begin{tikzpicture}[thick]
% 	\def\dx{0.2}
% 	\def\w{13.5}
% 	\def\dw{3.5}
% 	\def\h{3}
% \draw[rounded corners=10] (0,0) rectangle (\w,\h);
% \draw[rounded corners=10] (\dx,\dx) rectangle (\w-\dw,\h-\dx);
% \draw[rounded corners=10] (2*\dx,2*\dx) rectangle (\w-2*\dw,\h-2*\dx);
% \node[anchor=north east] at (\w,\h-3*\dx) {\parbox{3cm}{$\mathcal{C}_{m,L}$\\sector-bounded nonlinearity}};
% \node[anchor=north east] at (\w-\dw,\h-3*\dx) {\parbox{3cm}{$\mathcal{M}_{m,L}$\\slope-restricted nonlinearity}};
% \node[anchor=north east] at (\w-2*\dw,\h-3*\dx) {\parbox{5.4cm}{$\mathcal{F}_{m,L}$\\gradient of a strongly convex\\ func. with Lipschitz gradients}};
% \end{tikzpicture}
% \begin{caption}Venn diagram depicting nested classes of nonlinearities (or oracles).\label{fig:venn_diagram}
% \end{caption}
% \end{figure}
%%%%%%%%%%%%%%%%%%%%%%%%%%%%%%%%%%%%%%%%%%%%%%%%%%
\section{Certifying convergence rates using dissipativity}
\subsection{Certifying geometric convergence}
Geometric convergence, also known as \emph{exponential stability} in the controls literature and \emph{linear convergence} in the optimization literature, can be verified using a modified dissipation inequality of the form 
\begin{equation}\label{eq:diss_basic}
V(\xi^{k+1}) - q V(\xi^k) \leq S(\xi^k, u^k),
\end{equation}
where $0 \leq q < 1$. We distinguish between two important cases for use in algorithm analysis.
\begin{enumerate}
	\item If the supply rate is nonpositive, $S(\xi^k,u^k) \leq 0$, then the dissipation inequality implies that $V(\xi^{k+1}) \leq q V(\xi^k)$. So, $V$ decreases by a factor of $q$ at every timestep, which means that $V$ is a Lyapunov function that certifies geometric convergence. This case will be used when $\grad f \in \mathcal{C}_{m,L}$ or $\grad f \in \mathcal{M}_{m,L}$.
	\item If the supply rate satisfies the more general condition $S(\xi^k,u^k) \leq q \psi(\xi^k) - \psi(\xi^{k+1})$, where $\psi$ is a nonnegative function, then \eqref{eq:diss_basic} can be rewritten as
$V(\xi^{k+1}) + \psi(\xi^{k+1}) \leq q \left( V(\xi^k) + \psi(\xi^k) \right)$.
In other words, $V(x) + \psi(x)$ is a Lyapunov function that certifies geometric convergence. This case will be used when $\grad f \in \mathcal{F}_{m,L}$.
\end{enumerate}
\subsection{Certifying other convergence rates}
Dissipativity can also be used to certify other rates of convergence. For example, consider gradient descent (the algorithm state is $\xi^k \defeq x^k$), where the function value $f(x^k)$ decreases at each iteration. If the standard dissipation inequality~\eqref{eq:diss_will} holds with a supply rate that satisfies $S(x^k,u^k) \leq -\bigl(f(x^k) - f(x^\star)\bigr)$, then the dissipation inequality implies that $V(x^{k+1}) - V(x^k) \leq -\bigl(f(x^k) - f(x^\star)\bigr)$. Summing over $k$, 
\[
	f(x^k) - f(x^\star) \leq \frac{1}{k+1} \sum_{i=0}^{k} \left( f(x^i) - f(x^\star) \right)  \leq \frac{1}{k+1} \left( V(x^0) - V(x^{k+1}) \right) \leq \frac{1}{k+1} V(x^0).
\]
In other words, the algorithm converges in \emph{function value}, and the function value decreases at the \emph{sublinear rate} $1/k$, which is slower than geometric convergence.
Dissipativity can also be used to certify other sublinear rates, such as $1/k^2$ \cite{hu2017dissipativity}.
Ultimately, dissipation inequalities lead to Lyapunov functions. However, dissipation inequalities also provide a framework by which to treat uncertain disturbance inputs. Roughly, this is done by using supply rates that are satisfied by the disturbances and then searching for storage functions such that a dissipation inequality is satisfied (thereby certifying robust stability).
Even if the system in question is not mechanical in nature or there is no clear notion of \emph{energy}, the system may still satisfy a dissipation inequality. This is akin to the idea of that Lyapunov functions can be used to certify the stability of an autonomous differential or difference equation, even when the equation does not describe a physical system.
In the context of algorithm analysis, we view the outputs of the oracles as disturbance inputs for the iterative algorithm in question and use dissipativity theory to certify robust convergence properties.
%%%%%%%%%%%%%%%%%%%%%%%%%%%%%%%%%%%%%%%%%%%%%%%%
\section{Case study: Gradient Descent}
Gradient descent is perhaps the most recognizable iterative algorithm. We will now illustrate different approaches for analyzing gradient descent for a simple class of functions. Consider the (unconstrained) problem of minimizing $f(x)$, where $f$ is $m$-strongly convex and has $L$-Lipschitz gradients. That is, $f \in \mathcal{F}_{m,L}$. Gradient descent is defined by the iteration
\begin{equation}\label{eq:gdescent}
	x^{k+1} = x^k - \eta \grad f(x^k)
	\qquad\text{for }k=0,1,\dots.
\end{equation}
Our task is to find the worst-case linear convergence rate $\rho(\mathcal{A},\mathcal{F})$ and ultimately the choice of stepsize $\eta$ that leads to the fastest convergence rate.
\paragraph{Nesterov's analysis}
A classical approach to analyzing gradient descent is due to Nesterov~\cite[Thm.~2.1.15]{nesterov} and begins by bounding the error at the $(k+1)^\text{st}$ iterate in terms of the error at the $k^\text{th}$ iterate:
\begin{align}
\norm{x^{k+1} - x^\star}^2
&\overset{\eqref{eq:gdescent}}{=} \norm{x^k - x^\star - \eta \grad f(x^k)}^2 \notag \\
&= \norm{x^k - x^\star}^2 + \eta^2 \norm{\grad f(x^k)}^2 - 2\eta \left(x^k-x^\star\right)^\tp \grad f(x^k) \notag \\
&\overset{\eqref{eq:sector}}{\leq} \norm{x^k - x^\star}^2 + \eta^2 \norm{\grad f(x^k)}^2
-\tfrac{2\eta}{L+m}\left( m L \norm{x^k-x^\star}^2 + \norm{\grad f(x^k)}^2 \right) \notag\\
&= \left( 1 - \tfrac{2\eta mL}{L+m} \right) \norm{x^k - x^\star}^2 + \eta \left( \eta - \tfrac{2}{L+m} \right) \norm{\grad f(x^k)}^2. \label{eq:turtle}
\end{align}
If it is further assumed that $0 \leq \eta \leq \tfrac{2}{L+m}$, the second term in~\eqref{eq:turtle} is nonpositive, and it can be concluded that the error shrinks at every iteration according to
\[
\norm{x^{k+1} - x^\star} \leq \sqrt{ 1 - \tfrac{2\eta mL}{L+m} }\, \norm{x^k - x^\star}.
\]
The contraction factor $\sqrt{ 1 - \tfrac{2\eta mL}{L+m} }$ is an upper bound on the worst-case convergence rate. If the stepsize $\eta \in [0, \tfrac{2}{L+m}]$ is selected to minimize this upper bound, the optimal stepsize is $\eta = \frac{2}{L+m}$ and yields the bound $\rho \geq \tfrac{L-m}{L+m}$ on the worst-case convergence rate.
\paragraph{Polyak's analysis}
Another approach to analyzing gradient descent is due to Polyak~\cite[\S1.4, Thm.~3]{polyak}. Assume that $f$ is twice differentiable. By the fundamental theorem of calculus, 
\begin{align*}
\grad f(x^k) = \grad f(x^\star) + \int_0^1 \nabla^2f( x^\star + \tau(x^k-x^\star))(x^k-x^\star)\,\mathrm{d}\tau
= A_k (x^k - x^\star),
\end{align*}
where $A_k \defeq \int_0^1 \nabla^2 f(x^\star+\tau(x^k-x^\star))\,\mathrm{d}\tau$. Then, we can bound the error at the $(k+1)^\text{st}$ iterate in terms of the error at the $k^\text{th}$ iterate using the triangle inequality:
\[
	\norm{x^{k+1} - x^\star}
	= \norm{x^k - x^\star - \eta \grad f(x^k)}
	= \norm{ (I-\eta A_k) (x^k - x^\star)}
	\leq \norm{I-\eta A_k} \cdot \norm{x^k-x^\star}.
\]
It follows from~\eqref{eq:monotone3} that $m I \preceq \nabla^2 f(x) \preceq L I$ for all $x$. Therefore, $m I \preceq A_k \preceq L I$. Since $I-\eta A_k$ is symmetric, its norm can be bounded in terms of the smallest and largest eigenvalues of $A_k$, which gives a bound on the worst-case convergence rate:
\[
	\rho \geq \norm{I-\eta A_k} = \max\left\{ |1-\eta m|, |1-\eta L| \right\}.
	\]
This bound can be minimized if $\eta = \tfrac{2}{L+m}$, which yields $\rho \geq \tfrac{L-m}{L+m}$.
Both the Nesterov and Polyak analyses of gradient descent yield the same optimal stepsize of $\eta = \tfrac{2}{L+m}$ and the same worst-case convergence rate bound $\rho \geq \tfrac{L-m}{L+m}$. However, the bounds on $\rho$ disagree when $\eta \neq \tfrac{2}{L+m}$; Nesterov's bound is more conservative than Polyak's. Conversely, Nesterov's approach only assumed $\grad f$ was sector-bounded, while Polyak made the stronger assumptions that $f$ is twice differentiable and $\grad f$ is slope-restricted.
It is possible to prove Polyak's bound without assuming twice-differentiability, but it requires a different approach.
\paragraph{Dissipativity approach}
The gradient method can also be analyzed using dissipativity theory. First write the gradient method as a linear time-invariant dynamical system in feedback with $\grad f$:
\begin{subequations}\label{eq:graddyn}
\begin{align}
	x^{k+1} &= x^k - \eta u^k, \\
	y^k &= x^k, \\
	u^k &= \grad f( y^k ).
\end{align}
\end{subequations}
We then use the supply rate $S_{\mathcal{C}}$ associated with sector-bounded nonlinearities, given by \eqref{eq:sector}. The sector bound will not hold as written in~\eqref{eq:sector} because the oracle $u = \grad f(y)$ is zero when $y=x^\star = y^\star$ (the solution of the optimization problem), not when $y=0$. To account for this, use the inequality $S_\mathcal{C}(y-y^\star,u) \leq 0$. Although $y^\star$ may not be known in advance, the dissipativity approach only requires assuming existence of $y^\star$ and does not depend on its actual value.
The idea is to use a quadratic Lyapunov function candidate $V(x) = \norm{x - x^\star}^2$ and to certify a dissipation inequality of the form
\begin{equation}\label{eq:sectorlyap}
V(x^{k+1}) - \rho^2 V(x^k) \leq \lambda S_\mathcal{C}(y^k-y^\star,u^k),
\end{equation}
Where $\lambda \geq 0$. If~\eqref{eq:sectorlyap} holds, then $S_\mathcal{C}(y^k-y^\star,u^k) \leq 0$ implies that $V(x^{k+1}) \leq \rho^2 V(x^k)$, and therefore $\norm{x^{k+1} - x^\star} \leq \rho \norm{x^k-x^\star}$ and $\rho$ is an upper bound on the worst-case convergence rate. Substituting the definitions for $V$ and $S_\mathcal{C}$ and the dynamics~\eqref{eq:graddyn} into~\eqref{eq:sectorlyap}, 
\[
	\norm{x^k-x^\star-u^k}^2 - \rho^2 \norm{x^k-x^\star} \leq \lambda (m(x^k-x^\star) - u^k)^\tp (L(x^k-x^\star)-u^k).
\]
This quadratic expression in $x^k-x^\star$ and $u^k$ must hold for all choices of $x^k-x^\star$ and $u^k$, which implies the inequalities
\begin{equation}\label{eq:gradsdp1}
	\bmat{ 1-\rho^2 - \lambda m L & -\eta + \lambda \tfrac{L+m}{2} \\ -\eta + \lambda \tfrac{L+m}{2} & \eta^2 - \lambda } \preceq 0,\qquad \lambda \geq 0,
\end{equation}
where ``$\preceq$'' denotes inequality in the semidefinite sense. The inequalities~\eqref{eq:gradsdp1} do not depend explicitly on $y^\star$, even though the existence of $y^\star$ is assumed as part of the derivation.
The task of minimizing $\rho$ subject to~\eqref{eq:gradsdp1} is a semidefinite program (SDP). While typically solved numerically, SDPs can often be solved analytically in cases such as this one, where the matrices involved are small. In this case, a Schur complement of~\eqref{eq:gradsdp1} is used to obtain the equivalent inequalities:
\begin{equation}\label{eq:gradsdp1b}
	\rho^2 \geq 1- \lambda m L + \frac{(\eta - \lambda \tfrac{L+m}{2})^2}{\lambda - \eta^2},\qquad \lambda \geq \eta^2.
\end{equation}
Extremizing~\eqref{eq:gradsdp1b} with respect to $\lambda$ yields $\rho \geq \max\left\{ |1-\eta m|, |1-\eta L| \right\}$, which is the same bound as in Polyak's analysis.
For the simple case of gradient descent, the dissipativity approach produces the tightest possible bound (same as Polyak's bound), yet it only assumes $\grad f$ is sector-bounded. A nice feature of the dissipativity approach is that it is \emph{systematic}. When analyzing increasingly complicated algorithms, it becomes increasingly difficult to obtain useful convergence rate bounds via the traditional approach of ad hoc equation manipulation. The dissipativity approach provides a principled and scalable way to analyze algorithms.
%%%%%%%%%%%%%%%%%%%%%%%%%%%%%%%%%%%%%%%%%%%%%%%%
\section{Case study: Nesterov's Accelerated Method}
\begin{figure}[ht]
	\centering
	% \printinunitsof{in}\prntlen{\linewidth}
	\includegraphics{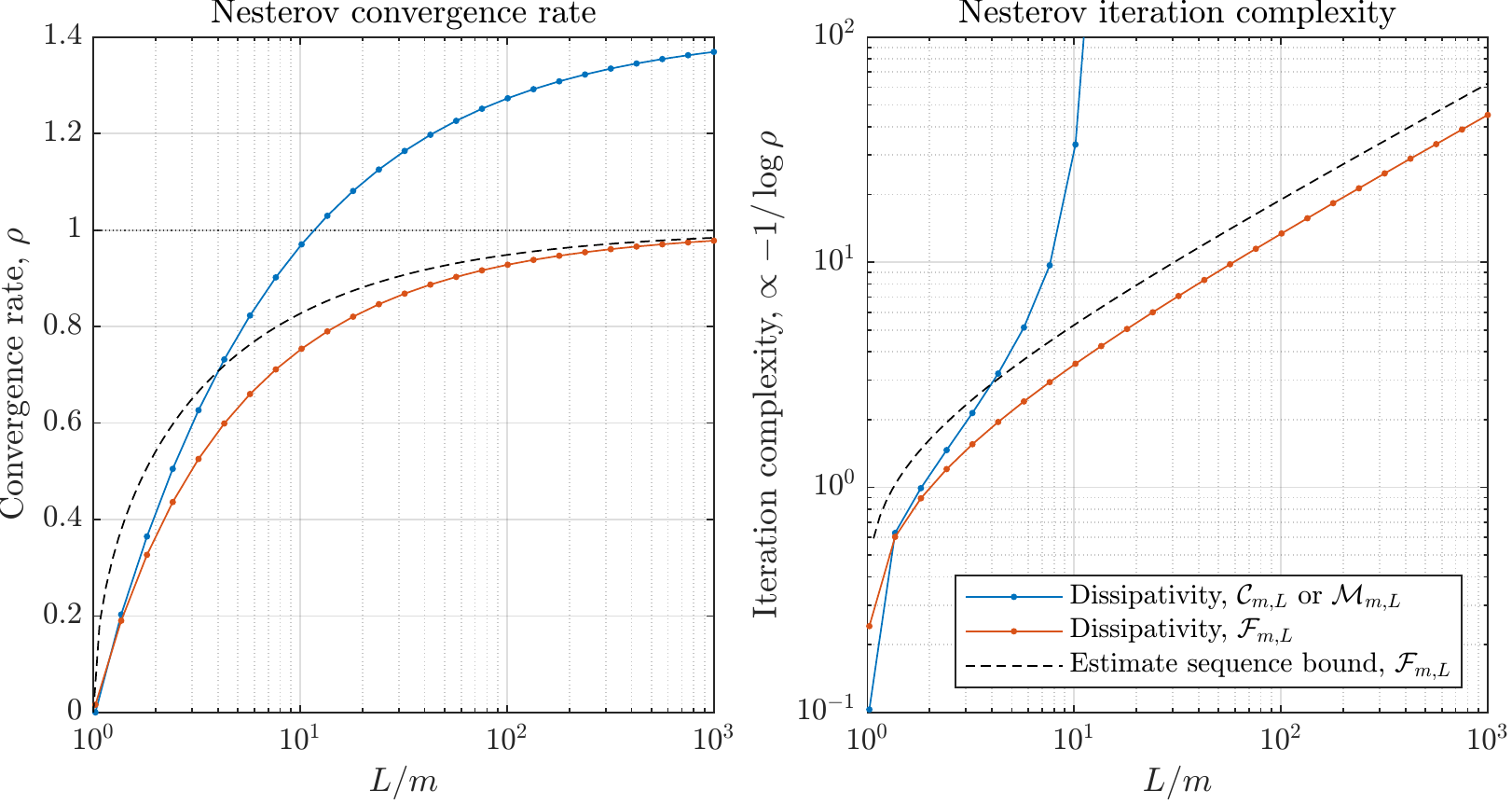}
	\caption{Worst-case convergence analysis of Nesterov's accelerated method with standard tuning. \textbf{Left:} Worst-case convergence rate $\rho$ as a function of the condition number $L/m$. In the case where $f$ is $m$-strongly convex with $L$-Lipschitz gradient ($\mathcal{F}_{m,L}$), the dissipativity approach yields an improvement on the bound found using estimate sequences~\cite[\S2.2.1]{nesterov}. Also shown is the case where $\grad f$ is only assumed to be sector-bounded ($\mathcal{C}_{m,L}$) or slope-restricted ($\mathcal{M}_{m,L}$), obtained using the dissipativity approach. Both cases have the same worst-case bound and are only guaranteed to converge ($\rho < 1$) when $L/m$ is relatively small. \textbf{Right:} The same data as the left plot, but instead we plot $-1/\log \rho$, which is proportional to the iteration complexity (the number of iterations required to ensure convergence to within a prespecified tolerance). As $\rho\to 1$ on the left plot (slower convergence), the iteration complexity tends to $+\infty$ on the right plot. When $\rho>1$, the algorithm may not converge, so iteration complexity is infinite. The dissipativity approach improves upon the estimate sequence bound by a constant factor of approximately $1.38$.}
	\label{fig:nestplot}
\end{figure}

\begin{figure}[ht]
	\centering
	% \printinunitsof{in}\prntlen{\linewidth}
	\includegraphics{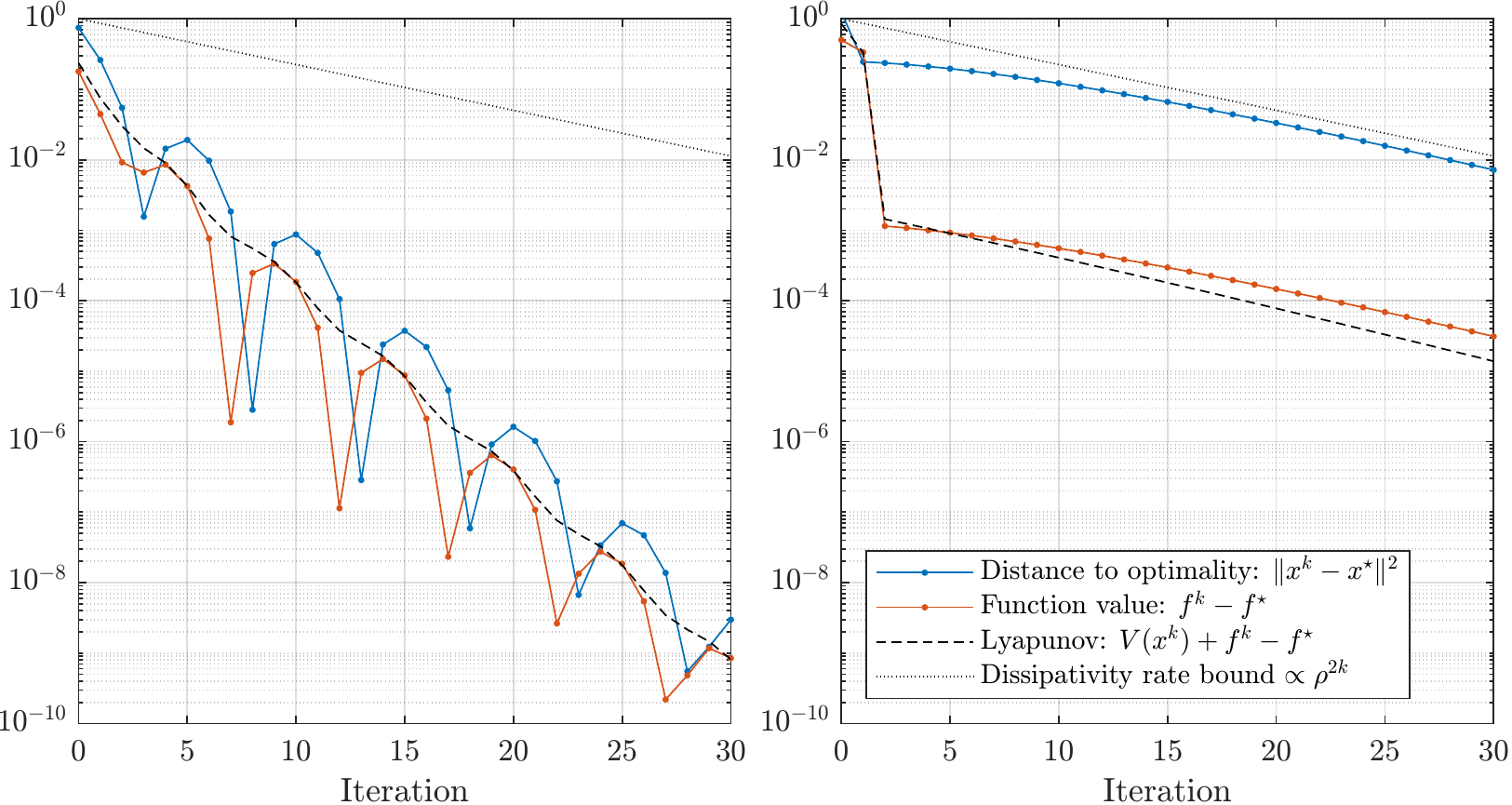}
	\caption{Simulation of Nesterov's accelerated method applied to two test functions that are $m$-strongly convex with $L$-Lipschitz gradients. In both cases, $L=1$ and $m=\tfrac{1}{100}$ were used and the squared distance to optimality $\norm{x^k-x^\star}^2$, the function value $f^k-f^\star$, and the Lyapunov function found using the dissipativity approach were plotted. \textbf{Left:} The function given in~\eqref{eq:f1}. In this case, neither the distance to optimality nor the function value decrease monotonically. However, the Lyapunov function does. Convergence is faster than that predicted by the dissipativity approach.
	\textbf{Right:} The function given in~\eqref{eq:f2}. In this case, the convergence is slower and the numerical bound from the dissipativity approach appears to be tight.
	For both functions, Nesterov's method was initialized with $x_1=1$ and $x_2=0.5$.}
	\label{fig:simplot}
\end{figure}

Consider Nesterov's accelerated method applied to a continuously differentiable function $f$, with access to the oracle $\grad f$. As in the gradient descent example, we will consider $f$ that is $m$-strongly convex with $L$-Lipschitz gradients. That is, $f\in \mathcal{F}_{m,L}$. It will be shown that the dissipativity approach can find the tightest known convergence rate for this algorithm. Recall that Nesterov's accelerated method is characterized by the iteration
\begin{subequations}\label{eq:nesterov_dynamics}
\begin{align}
	x_1^{k+1} &= (1+\beta)x_1^k - \beta x_2^k - \eta u^k, \\
	x_2^{k+1} &= x_1^k, \\
	y^k &= (1+\beta)x_1^k - \beta x_2^k.
\end{align}
\end{subequations}
The classical approach to analyzing Nesterov's accelerated method is \emph{estimate sequences}~\cite[\S2.2.1]{nesterov}, which are a recursively generated sequence of quadratic bounds on the worst-case convergence rate $\rho$. The approach is rather involved, so the exposition is omitted. The net result is that the worst-case convergence rate satisfies the bound
\begin{equation}\label{eq:nest_bound_estimate_sequence}
	\rho \geq \sqrt{1 - \sqrt{\tfrac{m}{L}}}.
\end{equation}
We showcase the versatility of the  dissipativity approach by analyzing the cases where $\grad f$ belongs to $\mathcal{C}_{m,L}$, $\mathcal{M}_{m,L}$, or $\mathcal{F}_{m,L}$. In all of these cases, as in the analysis of gradient descent, we use a change of variables to shift the optimal point to zero. That is, the dynamics are re-expressed in terms of $x^k-x^\star$ and $y^k-y^\star$. Since the convergence rate bound found via dissipativity is independent of the optimal point, the notation is simplified by assuming without loss of generality that the optimal point is at zero, so $x^\star = y^\star = u^\star = 0$.
\paragraph{Nesterov for sector-bounded gradients}
For the sector-bounded case, the same approach is used as in the analysis of gradient descent. There are two states, $x_1$ and $x_2$, so the storage function will be a positive-definite quadratic function of both:
\[
V(x) = x^\tp P x,
\quad\text{where }
x = \bmat{x_1 \\ x_2}\text{ and }P \succ 0.
\]
The supply rate $S_\mathcal{C}$ defined in~\eqref{eq:sector} is used, and we seek the smallest $\rho$ such that the following dissipation inequality is satisfied for some $\lambda \geq 0$:
\begin{equation}\label{eq:diss_nest}
V(x^{k+1}) - \rho^2 V(x^k) \leq \lambda S_\mathcal{C}(y^k,u^k).
\end{equation}
Since the presence of both $P$ and $\lambda$ renders this dissipation inequality homogeneous, it can be normalized by setting $\lambda=1$. Upon substituting the dynamics~\eqref{eq:nesterov_dynamics}, the dissipation inequality~\eqref{eq:diss_nest} becomes a quadratic inequality in $(x^k, u^k)$, which leads to a semidefinite program (as in the case study for gradient descent).
\paragraph{Nesterov for slope-restricted gradients}
For the slope-restricted case, the aim is to use the supply rate $S_\mathcal{M}$ defined in~\eqref{eq:monotone}. However, using consecutive iterates as the two points in the storage function [for example $(y^k,u^k)$ and $(y^{k+1},u^{k+1})$] will cause the dissipation inequality to depend on both $u^k$ and $u^{k+1}$. Therefore, the storage function must be augmented to depend on both $x^k$ and $x^{k+1}$. For the supply rate, there are many choices. In particular, apply~\eqref{eq:monotone} at any pair of points chosen among $\{y^k, y^{k+1},y^\star\}$. Since the expression for $S_\mathcal{M}$ is symmetric in its arguments, this leads to three supply rate inequalities:
\begin{gather*}
	S_\mathcal{M}^1 \defeq S_\mathcal{M}( y^k, y^{k+1}, u^k, u^{k+1}) \leq 0, \qquad
	S_\mathcal{M}^2 \defeq S_\mathcal{M}( y^k, y^\star, u^k, u^\star) \leq 0,\\
	S_\mathcal{M}^3 \defeq S_\mathcal{M}( y^{k+1},y^\star, u^{k+1},u^\star) \leq 0.
\end{gather*}
Ultimately, the dissipation inequality is
\begin{equation}\label{eq:diss_nest2}
	V(x^{k+1},x^{k+2}) - \rho^2 V(x^{k},x^{k+1})
	\leq \sum_{i=1}^3 \lambda_i S_\mathcal{M}^i,
\end{equation}
where $\lambda_1,\lambda_2,\lambda_3 \geq 0$ and $V(\cdot, \cdot)$ is a positive-definite quadratic. Upon substituting the dynamics~\eqref{eq:nesterov_dynamics}, the dissipation inequality~\eqref{eq:diss_nest2} becomes a quadratic inequality in $(x^k,u^k,u^{k+1})$, which again leads to a semidefinite program.
Further augmenting the storage function to include more consecutive iterates $\{x^k,x^{k+1},\dots,x^{k+r}\}$ will further increase the number of supply rate inequalities available for use, potentially yielding less conservative upper bounds on the worst-case convergence rate $\rho$.
\paragraph{Nesterov for gradients of strongly convex functions}
For the case $\grad f \in \mathcal{F}_{m,L}$, the aim is to use the supply rate $S_{\mathcal{F}}$ defined in~\eqref{eq:strcvx}. For brevity, we write $f^k \defeq f(y^k)$ and $f^\star \defeq f(y^\star)$. This time, the supply rate is not symmetric in its arguments, so there are six inequalities to choose from:
\begin{align*}
	S_\mathcal{F}^1 &\defeq S_\mathcal{F}( y^{k+1}, y^{k}, u^{k+1}, u^{k}) \leq f^k - f^{k+1}, &
	S_\mathcal{F}^2 &\defeq S_\mathcal{F}( y^k, y^{k+1}, u^k, u^{k+1}) \leq f^{k+1} - f^k,\\
	S_\mathcal{F}^3 &\defeq S_\mathcal{F}( y^k, y^\star, u^k, u^\star) \leq f^\star - f^k, &
	S_\mathcal{F}^4 &\defeq S_\mathcal{F}( y^\star, y^k, u^\star, u^k) \leq f^k - f^\star, \\
	S_\mathcal{F}^5 &\defeq S_\mathcal{F}( y^{k+1},y^\star, u^{k+1},u^\star) \leq f^\star - f^{k+1}, &
	S_\mathcal{F}^6 &\defeq S_\mathcal{F}( y^\star,y^{k+1}, u^\star,u^{k+1}) \leq f^{k+1} - f^\star.
\end{align*}
In this context, the storage function $V$ will not itself be a Lyapunov function. Rather, the Lyapunov function will be of the form $V(\cdot) + f^k$. Therefore, the storage function $V$ need not be positive definite. We therefore use two inequalities:
\begin{subequations}\label{eq:diss_nest3}
\begin{align}\label{eq:diss_nest3a}
	V(x^{k+1},x^{k+2}) - \rho^2 V(x^k,x^{k+1})
	&\leq \sum_{i=1}^6 \lambda_i S_\mathcal{F}^i,\\
	-V(x^{k},x^{k+1}) &\leq \sum_{i=1}^6 \mu_i S_\mathcal{F}^i,
	\label{eq:diss_nest3b}
\end{align}
\end{subequations}
where the $\lambda_i$ and $\mu_i$ are nonnegative. Ultimately, the aim is for the supply rates to satisfy
\begin{subequations}\label{eq:falcon}
\begin{align}
\sum_{i=1}^6 \lambda_i S_\mathcal{F}^i &\leq \rho^2 \left( f^k - f^\star \right) - \left( f^{k+1}-f^\star\right), \\
\sum_{i=1}^6 \mu_i S_\mathcal{F}^i &\leq \left( f^k - f^\star \right),
\end{align}
\end{subequations}
which is ensured via the additional linear constraints
\begin{align*}
	\lambda_1 - \lambda_2 - \lambda_3 + \lambda_4 &= \rho^2, &
	-\lambda_1 + \lambda_2 - \lambda_5 + \lambda_6 &= -1, \\
	\mu_1 - \mu_2 - \mu_3 + \mu_4 &= 1, &
	-\mu_1 + \mu_2 - \mu_5 + \mu_6 &= 0.
\end{align*}
Combining~\eqref{eq:diss_nest3} with \eqref{eq:falcon}, it follows that $V(x^{k},x^{k+1})+(f^k-f^\star)$ is a Lyapunov function that certifies geometric convergence with rate $\rho$. As with the case $\mathcal{M}_{m,L}$, it is possible to further augment the storage function and use more supply rates, which can potentially yield less conservative upper bounds on the worst-case convergence rate $\rho$.
\paragraph{Numerical simulation}
When examining the semidefinite programs~\eqref{eq:diss_nest}, \eqref{eq:diss_nest2}, and \eqref{eq:diss_nest3}, they differ from the gradient descent case~\eqref{eq:gradsdp1} in that they are not linear in $\rho^2$ due to the presence of the product $\rho^2 P$. However, they are linear in $P$ and the  $\lambda_i$'s and $\mu_i$'s for each fixed $\rho$, so they can be solved by bisection on $\rho$.
Implementing these solutions with the default tuning of Nesterov's method, which is $\eta = \frac{1}{L}$ and $\beta = \frac{\sqrt{L} - \sqrt{m} }{\sqrt{L}+\sqrt{m}}$, we obtain the results displayed in Figure~\ref{fig:nestplot}. These results mirror those reported in~\cite{lessard2016analysis}, although that article used the theory of integral quadratic constraints (IQCs)~\cite{megrantzer} and Zames--Falb multipliers~\cite{HeaWilZF,ZamFal} instead of dissipativity. IQC theory is closely related to dissipativity theory~\cite{seiler_diss}, and the dissipativity approach presented above is algebraically equivalent to the approach used in~\cite{lessard2016analysis}.
As shown in Figure~\ref{fig:nestplot}, the worst-case rate certified by the dissipativity approach improves upon the rate found using the classical approach of estimate sequences. Figure~\ref{fig:nestplot} also plots the \emph{iteration complexity}, which is the number of iterations required to reach a certain error $\varepsilon$. The iteration complexity is proportional to $-1/\log\rho$. The results in Figure~\ref{fig:nestplot} were obtained numerically using CVX~\cite{cvx}, a Matlab package for specifying and solving convex programs. In this case, the semidefinite program is more complicated than the one for gradient descent, yet still simple enough to be solved analytically~\cite{nesterov_analytic}.
When using the Lyapunov function $V(x^{k+1},x^k) + (f^k-f^\star)$, it is not required that $V$ itself be positive, nor is it required that $V$ or $f^k$ be monotonically decreasing. Figure~\ref{fig:simplot} shows two examples of Nesterov's method applied to functions in $\R^2$. Both cases use functions with $\grad f \in \mathcal{F}_{m,L}$ with the same values of $(m,L)$ and the same initialization. The functions used are
\begin{subequations}
	\begin{align}
		\label{eq:f1}
	f_1(x_1,x_2) &\defeq \frac{m}{2}(x_1^2+x_2^2) + (L-m) \log\left( e^{-x_1} + e^{x_1/3+x_2} + e^{x_1/3-x_2}\right) \\
	f_2(x_1,x_2) &\defeq \tfrac{L}{2}x_1^2 + \tfrac{m}{2}x_2^2.
\label{eq:f2}
	\end{align}
\end{subequations}
In the first case (left panel of Figure~\ref{fig:simplot}),  the distance to optimality $\norm{x^k-x^\star}$ and the function value $f(x^k) - f^\star$ both converge nonmonotonically. The Lyapunov function found using the dissipativity approach is, however, monotone. The dissipativity interpretation is that some of the energy is stored in the function value, while some is stored in the states. Energy sloshes back and forth between both. However, the total energy is dissipated at a rate bounded by $\rho^{2k}$. In this case, the bound is loose; the algorithm converges significantly faster than its worst-case bound. In the second case (right panel of Figure~\ref{fig:simplot}), the energy is dissipated more slowly and matches the worst-case bound.
%%%%%%%%%%%%%%%%%%%%%%%%%%%%%%%%%%%%%%%%%%%%%%%%
\section{Case study: Alternating Direction Method of Multipliers}
\begin{figure}[ht]
	\centering
	\includegraphics{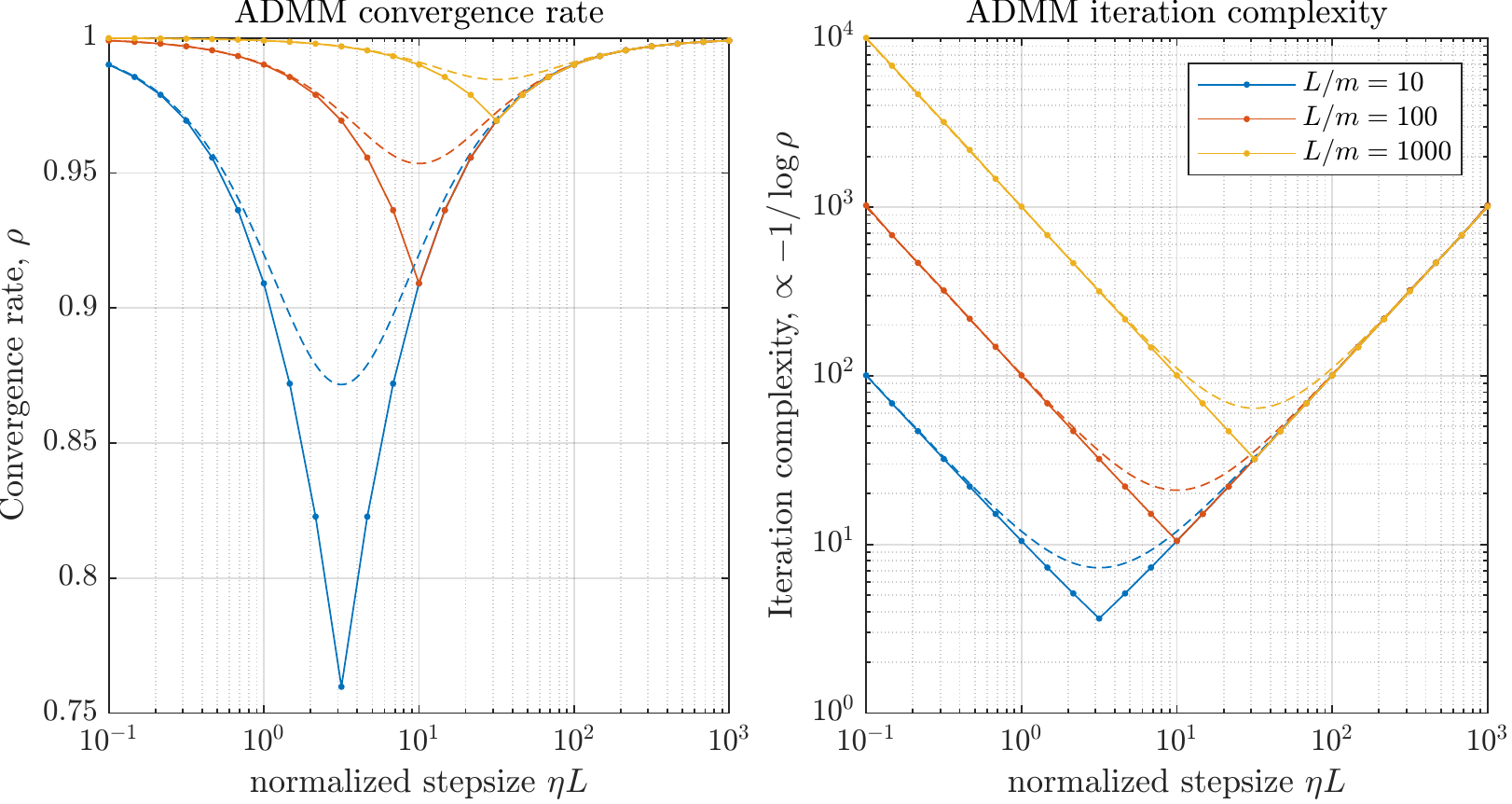}
	\caption{Worst-case convergence analysis of Alternating Direction Method of Multipliers (ADMM) applied to the problem of minimizing $f(x)+g(x)$, where $f$ is $m$-strongly convex with $L$-Lipschitz gradient and $g$ is convex but not necessarily differentiable. This includes the case of constrained optimization (where $g$ is the indicator function of a convex set). Solid lines show the upper bound found by solving the dissipation inequality~\eqref{eq:admm_sdp}. Dashed lines show the upper bound~\eqref{yinub} from~\cite[Cor.~3.1]{deng2016global}, which is looser than the SDP bound.
	\textbf{Left:} Worst-case convergence rate $\rho$ as a function of the normalized stepsize $\eta L$ for different choices of $L/m$. \textbf{Right:} The same data as the left plot, but instead plot $-1/\log \rho$, which is proportional to the iteration complexity (the number of iterations required to ensure convergence to within a prespecified tolerance).}\label{fig:admmplot}
\end{figure}

Consider the ADMM algorithm applied to a composite optimization problem $\min_x f(x) + g(x)$, as described by the interconnections shown in Figure~\ref{fig:blkdiag_admm}. We investigate how to tune the stepsize $\eta$ to ensure the fastest possible worst-case convergence.
Assume $f$ is strongly convex and $g$ is convex. In other words, $\grad f \in \mathcal{F}_{m,L}$ and $\partial g \in \mathcal{M}_{0,\infty}$. We can obtain an upper bound for the worst-case convergence rate using the supply rates for $\mathcal{C}_{m,L}$ and $\mathcal{C}_{0,\infty}$, respectively. Using the standard Lyapunov candidate
\[
V(x) = x^\tp P x,
\quad\text{where }
x \defeq \bmat{z \\ w}\text{ and }P \succ 0,
\]
and a similar dissipation inequality to~\eqref{eq:sectorlyap}, which is used to analyze the gradient method:
\begin{equation}\label{eq:sectorlyap2}
	V(x^{k+1}) - \rho^2 V(x^k) \leq \lambda_1 S_\mathcal{C}^{m,L}(y_1^k-y_1^\star,u_1^k) + \lambda_2 S_\mathcal{C}^{0,\infty}(y_2^k-y_2^\star,u_2^k).
\end{equation}
In~\eqref{eq:sectorlyap2}, we used superscripts with $S_\mathcal{C}$ to denote the bounds of the sector. For the case $S_\mathcal{C}^{0,\infty}$, which corresponds to the set $\mathcal{C}_{0,\infty}$, we divide~\eqref{eq:sector} by $L$ and take the limit $L\to \infty$, which yields
\[
S_\mathcal{C}^{0,\infty}(y,u) = \bmat{y \\ u}^\tp \bmat{ 0 & -\tfrac{1}{2} \\ -\tfrac{1}{2} & 0 } \bmat{y \\ u }.	
\]
Upon substituting the dynamics~\eqref{eq:admm_grad} into the dissipation inequality~\eqref{eq:sectorlyap2}, a quadratic inequality is obtained in the variables $\{z^k-z^\star, w^k-w^\star, u_1^k, u_2^k\}$, which yields the semidefinite inequality
\begin{multline}\label{eq:admm_sdp}
\bmat{\star}^\tp P \bmat{ 1 & 0 & -\eta & -\eta \\ 0 & 0 & 0 & \eta}
- \rho^2 \bmat{\star}^\tp P \bmat{ 1 & 0 & 0 & 0 \\ 0 & 1 & 0 & 0} \\
\preceq  \lambda_1 \bmat{\star}^\tp \bmat{ m L & -\tfrac{L+m}{2} \\ -\tfrac{L+m}{2} & 1} \bmat{ 1 & -1 & -\eta & 0 \\ 0 & 0 & 1 & 0} + \lambda_2 \bmat{\star}^\tp \bmat{ 0 & -\tfrac{1}{2} \\ -\tfrac{1}{2} & 0} \bmat{ 1 & 0 & -\eta & -\eta \\ 0 & 0 & 0 & 1},
\end{multline}
where $\star$ denotes the term required to make the quadratic forms symmetric. For example, in $\bmat{\star}^\tp P X$, $\star = X$.
For each fixed $\rho,\eta,m,L$, \eqref{eq:admm_sdp} is an SDP in the variables $P \succ 0$ and $\lambda_1, \lambda_2 \geq 0$. For different choices of stepsize $\eta$ and condition number $L/m$, \eqref{eq:admm_sdp} is solved using a bisection search to find the smallest $\rho$ that ensures feasibility. Numerically obtained SDP upper bounds are plotted in Figure~\ref{fig:admmplot}, alongside the analytical bound from~\cite[Cor.~3.1]{deng2016global}, which is given by \begin{equation}\label{yinub}
\rho_\textup{ub} \leq \sqrt{\frac{\eta ^2 L m+1}{\eta ^2 L m+2 \eta  m+1}}.	
\end{equation}
The convergence rate upper bounds shown in Figure~\ref{fig:admmplot} may be conservative because supply rates $\mathcal{C}_{m,L}$ and $\mathcal{C}_{0,\infty}$ were used rather than the more precise $\mathcal{F}_{m,L}$ and $\mathcal{M}_{0,\infty}$, respectively. Nevertheless, a matching lower bound can be found by considering a specific quadratic problem instance, as in~\cite{nishihara2015general,ghadimi2014optimal}. Consider the composite objective with $f(x) = \frac{1}{2}x^\tp Q x$, where the largest and smallest eigenvalues of $Q$ are given by $L$ and $m$, respectively, and $g(x) = \frac{1}{2}\delta \norm{x}^2$ with $\delta \geq 0$. In this case, ~\eqref{eq:admm_prox} reduces to the linear update
$
z^{k+1} = (1+\eta\delta)^{-1}(I+\eta Q)^{-1}(I+\eta^2\delta Q)	z^k
$.
If $\lambda \in [m,L]$ is an eigenvalue of $Q$, the eigenvalues of the update matrix for $z^k$ are $(1+\eta \delta)^{-1}(1+\eta \lambda)^{-1}(1+\eta^2\delta \lambda)$. Therefore, the greatest lower bound is found by extremizing:
\[
\rho_\textup{lb} \geq \sup_{\delta \geq 0} \sup_{m\leq \lambda\leq L} \frac{1+\eta^2\delta \lambda}{(1+\eta\delta)(1+\eta\lambda)}
\geq \max\left( \frac{1}{1+\eta m}, \frac{\eta L}{1+\eta L} \right),
\]
where the second inequality is found by picking $\delta=0$ and $\delta\to\infty$. This lower bound coincides precisely with the numerical SDP upper bound plotted in Figure~\ref{fig:admmplot}, which suggests that the upper bound was tight.
A rigorous proof that the upper and lower bounds match requires finding an analytic solution to \eqref{eq:admm_sdp}. This SDP is larger and has more variables than the gradient descent SDP~\eqref{eq:gradsdp1}, which is why analytic solutions are more difficult to obtain in this case. By inspection, an analytic solution for the peaks of interest in Figure~\ref{fig:admmplot} can be found, which are the stepsizes $\eta$ at which the worst-case convergence rate $\rho$ is fastest:
\begin{gather*}
	\eta = \frac{1}{\sqrt{L m}},\qquad \rho = \frac{\sqrt{L}}{\sqrt{L}+\sqrt{m}},\qquad
	\lambda_1 = 1, \qquad \lambda_2 = \frac{(L-m)^2}{L}, \\
P = \frac{\sqrt{m}(\sqrt{L} + \sqrt{m})(L-m)}{2}\bmat{1 & \sqrt{\frac{m}{L}} \\ \sqrt{\frac{m}{L}} & 1}.
\end{gather*}
The SDP analysis above was carried out using the gradient oracles~\eqref{eq:admm_grad}. The same worst-case rates are obtained if we instead used the proximal oracle formulation~\eqref{eq:admm_prox}, provided suitable adjustments are made to the supply rates. For example, monotonicity of the subgradient ($\partial g \in \mathcal{M}_{0,\infty}$) corresponds to firm nonexpansiveness of the proximal operator ($\prox_{\eta g} \in \mathcal{M}_{0,1}$).
Although the optimized convergence rate for ADMM, $\frac{\sqrt{L}}{\sqrt{L}+\sqrt{m}}$, is faster than that of gradient descent, $\frac{L-m}{L+m}$, the cost of a single iteration may be dramatically different in terms of wall-clock time, since computing one gradient  is likely far cheaper than computing one proximal operation. Indeed, as $\eta\to\infty$, $\prox_{\eta f}(x) \to \argmin_x f(x)$. So, when $\eta$ is large, a single proximal operation solves the unconstrained optimization problem.
%%%%%%%%%%%%%%%%%%%%%%%%%%%%%%%%%%%%%%%%%%%%%%%%
\section{General scalable algorithm analysis}
The approaches developed in the three previous case studies extend to algorithms in the general form~\eqref{eq:genalg}. In other words, assume an algorithm with updates of the form
\begin{align}\label{eq:gendyn}
\bmat{\xi^{k+1} \\ y_1^k \\ \vdots \\ y_m^k} &=
\bmat{ A & B_1 & \cdots & B_m \\ C_1 & D_{11} & \cdots & D_{1m} \\
\vdots & \vdots & \ddots & \vdots \\
C_m & D_{m1} & \cdots & D_{mm}}
\bmat{ \xi^k \\ u_1^k \\ \vdots \\ u_m^k}
\end{align}
and oracles $u_i^k = \phi_i(y_i^k)$ for $i=1,\dots,m$. Suppose we want to use a window of $r$ consecutive timesteps in our analysis. Then, all relevant supply rates for each of the oracles are written as in the case study on Nesterov's method. These will depend on the pairs $(y^\star,u^\star), (y^k,u^k),\dots,(y^{k+r},u^{k+r})$. These supply rates are called $\{S^1,\dots,S^M\}$. As before, choose a quadratic storage function $V(\xi^k,\dots,\xi^{k+r})$. To simplify the exposition, consider the case where the oracles are sector-bounded or slope-restricted. This leads to the  following inequality, similar to~\eqref{eq:diss_nest} and~\eqref{eq:diss_nest2}:
% \begin{subequations}\label{eq:diss_gen3}
	\begin{align}\label{eq:diss_gen3a}
		V(\xi^{k+1},\dots,\xi^{k+r+1}) - \rho^2 V(\xi^{k},\dots,\xi^{k+r})
		&\leq \sum_{i=1}^M \lambda_i S^i,
		% \\
		% -V(\xi^{k},\dots,\xi^{k+r}) &\leq \sum_{i=1}^M \mu_i S^i,
		% \label{eq:diss_gen3b}
	\end{align}
% \end{subequations}
where the $\lambda_i$ are nonnegative. If some of the oracles are gradients of strongly convex functions instead, we also include the positivity inequality as in~\eqref{eq:diss_nest3b} and the associated linear constraints on the $\lambda_i$ and $\mu_i$. The supply rates are quadratic functions of pairs of inputs and outputs, and upon substituting the dynamics~\eqref{eq:gendyn} to eliminate $\xi^{k+1},\dots,\xi^{k+r+1}$. The inequality~\eqref{eq:diss_gen3a} reduces to a quadratic expression in the variables $\{\xi^k,u^k,\dots,u^{k+r}\}$, where $u^{k+i} \defeq (u^{k+i}_1,\dots,u^{k+i}_m)$.
How large is the semidefinite program associated with~\eqref{eq:diss_gen3a}? If it is assumed that the algorithm has $n$ states and the oracles each map $\phi_i:\R^{d} \to \R^{d}$, then $A \in \R^{nd\times nd}$, $B_i \in \R^{nd\times d}$, $C_i \in \R^{d\times nd}$, and $D_{ij} \in \R^{d\times d}$. The storage function depends on $\xi^k,\dots, \xi^{k+1}$, so the quadratic form can be represented using a symmetric matrix $P\in\R^{n(r+1)d\times n(r+1)d}$. Note also the scalar variables $\lambda_1,\dots,\lambda_M$. Meanwhile, the entire inequality~\eqref{eq:diss_gen3a} depends on $\{\xi^k,u^k,\dots,u^{k+r}\}$, so it is a semidefinite constraint of size $(n+r+1)d \times (n+r+1)d$. 
\subsection{Scalability}
For typical iterative algorithms, $n$ is small; $n=1$ for gradient descent and $n=2$ for Nesterov's method. Tight convergence guarantees can also be obtained with relatively small $r$; our case studies required $r=1$ or $r=2$. However, in machine learning applications, it is not uncommon for $d$ to be very large (millions or billions), which would make the semidefinite program described above prohibitively large.
Although $d$ may be large, algorithms typically have highly structured system matrices $(A,B,C,D)$. For example, it is shown from~\eqref{eq:gdescent} that gradient descent has diagonal system matrices:
\begin{align*}
	A&=I_d, & B&=-\eta I_d, & C &= I_d, & D &= 0,
\end{align*}
where $I_d \in \R^{d\times d}$ is the identity matrix. Similarly, it is apparent from~\eqref{eq:nesterov_dynamics} that Nesterov's method has block-diagonal system matrices:
\begin{align*}
	A&=\bmat{1+\beta & -\beta \\ 1 & 0}\otimes I_d, &
	B&= \bmat{-\eta \\ 0} \otimes I_d, &
	C&= \bmat{1+\beta & -\beta} \otimes I_d, &
	D&= 0,
\end{align*}
where $\otimes$ denotes the Kronecker product. A similar structure is present in the supply rates. Therefore, the $P$ matrix from the storage function may also be written as $P = \hat P \otimes I_d$, where $\hat P \in \R^{n(r+1)\times n(r+1)}$. Thus, the semidefinite program decouples into $d$ identical (and much smaller) semidefinite programs, meaning that the dissipativity-based worst-case algorithm analysis only requires solving small semidefinite programs whose size depends on $n$ and $r$, but not $d$.
%%%%%%%%%%%%%%%%%%%%%%%%%%%%%%%%%%%%%%%%%%%%%%%%
\section{Concluding remarks}
The case studies above show how dissipativity can be applied to the analysis of iterative optimization algorithms. Further examples can be analyzed in an analogous manner, including distributed optimization algorithms~\cite{distralg,bcon_CDC}, stochastic and variance-reduction algorithms~\cite{sgdn, disstoch_ICML, hu2017unified}, and alternating algorithms for smooth games~\cite{giqc}.
The benefit of using dissipativity for algorithm analysis is that it provides a principled and modular framework where algorithms and oracles can be interchanged and analyzed depending on the case at hand. Additionally, the dissipativity approach is computationally tractable. The semidefinite programs solved in the case studies are small, with fewer than 20 variables. Importantly, the size of the semidefinite programs is independent of the dimension of the domain of the objective function.
\section{Acknowledgments}
This material is based upon work supported by the National
Science Foundation under Grants No. 1656951, 1750162, 1936648.
%%%%%%%%%%%%%%%%%%%%%%%%%%%%%%%%%%%%%%%%%%%%%%%%
\bibliographystyle{abbrv}
\begin{small}
\bibliography{algopt}
\end{small}

\appendix

\newpage
\section[50 years of dissipativity in algorithm analysis]{50 years of dissipativity in algorithm analysis}
\label{sb:history}
Lyapunov theory~\cite{lyapunov} is the main tool for the stability analysis of dynamical systems and has become a cornerstone of control theory as well.
Optimization algorithms have a natural interpretation as dynamical systems and therefore can also be analyzed using Lyapunov theory. The idea of viewing iterative algorithms as dynamical systems started perhaps in the early 1970s with the work of Tsypkin~\cite{tsypkin1971adaptation}. Tsypkin draws many parallels between optimization algorithms and control systems, with an emphasis on gradient-based algorithms. The main connections are outlined in the table below.\looseness=-1

\bigskip
\noindent
\begin{tabular}{p{7cm} p{1cm} p{7cm}}
	\toprule
	\textbf{Optimization} && \textbf{Controls} \\\midrule
	algorithm converges to a minimizer && equilibrium point is asymptotically stable \\[3mm]
	algorithm converges to a minimizer for all functions in a given class && robust stability \\[9mm]
	designing an algorithm with bounds on its worst-case performance && robust controller synthesis \\
	\bottomrule
\end{tabular}
\bigskip

At around the same time, Willems published his seminal work on dissipativity~\cite{willems1,willems2}, which generalized Lyapunov theory to systems with inputs, and he also generalized previous robust stability criteria such as passivity theory and the small-gain theorem~\cite{zames}. These ideas are at the core of nonlinear systems theory and endure to this day~\cite{khalil2002nonlinear}.
Another key body of work is Polyak's \emph{Introduction to Optimization}~\cite{polyak}, which covers the fundamentals of iterative methods for continuous and convex optimization. Polyak's book adopts a dynamical systems perspective, and even features a chapter on Lyapunov's method.
The problem of worst-case algorithm analysis is fundamentally a Lur'e problem~\cite{lure_postnikov}. An important distinction is that in the robust control literature, the goal is typically to prove \emph{stability}, whereas in algorithm analysis, the goal is to quantify the rate of convergence (which is akin to controlling the rate of exponential convergence). The modern tool for tackling this problem is integral quadratic constraints (IQCs), either in the frequency domain \cite{megrantzer,veenman2016robust} or in the time domain via dissipativity~\cite{seiler_diss}.
More recent works study algorithm analysis through the dynamical systems viewpoint have used IQC and dissipativity tools to achieve the tightest known bounds on many popular algorithms~\cite{lessard2016analysis,nishihara2015general,hu2017dissipativity}. These ideas have also been extended to synthesis, in an effort to design new optimization algorithms in a principled way~\cite{michalowsky2020robust,distralg,van2017fastest,cyrus2018robust}
Although this survey focused mostly on gradient-based iterative methods for continuous optimization, control perspectives have also been applied to (1) iterative methods for solving linear systems, (2) algorithms for solving ordinary differential equations, and (3) algorithms for solving linear programs. For a comprehensive survey of these topics, see~\cite{bhaya2006control}.
%%% SECOND SIDEBAR
\clearpage
\newpage
\section[Proving inequalities via the S-procedure]{Proving inequalities via the S-procedure}
\label{sb:sproc}
In the optimization literature, inequalities such as~\eqref{eq:sector2} are typically proved on a case-by-case basis using clever combinations of the triangle inequality, Cauchy--Schwarz, and other inequalities. Instead, the more fundamental property of positivity can often be used, namely if $S \leq 0$ is known to be true and $S = Q + P$ is split, where $P\geq 0$, then it follows that $Q \leq 0$.
Start with the definition
\[
	S_\mathcal{C}(y,u) \defeq \bmat{y\\u}^\tp\bmat{ mL & -\frac{L+m}{2} \\ -\frac{L+m}{2} & 1}\bmat{y\\u}
\]
and consider the algebraic identity
\begin{subequations}\label{eq:sproc_identity_all}
\begin{equation}
S_\mathcal{C}(y,u) = \tfrac{L-m}{2m} \bigl(  m^2 \norm{y}^2 - \norm{u}^2 \bigr) + \tfrac{L+m}{2m} \norm{my - u}^2,
\label{eq:sproc_identity}
\end{equation}
which can be directly verified to hold for all $y$ and $u$ by expanding both sides and comparing like terms. Since $0 < m \leq L$ by assumption, the second term in~\eqref{eq:sproc_identity} is nonnegative, as norms are always nonnegative. Therefore, by the positivity property, if $S_\mathcal{C}(y,u) \leq 0$, then $m^2 \norm{y}^2 -\norm{u}^2 \leq 0$. In other words, $m \norm{y} \leq \norm{u}$, which is the first half of~\eqref{eq:sector2}. The remaining inequalities in~\eqref{eq:sector2}--\eqref{eq:sector3} can be proven using a similar approach, with the further algebraic identities:
\begin{align}
S_\mathcal{C}(y,u) &= \tfrac{L-m}{2L} \bigl( \norm{u}^2 - L^2\norm{y}^2 \bigr) + \tfrac{L+m}{2L} \norm{u - Ly}^2, \\
S_\mathcal{C}(y,u) &= (L-m) \bigl(  m \norm{y}^2 - u^\tp y \bigr) + \norm{my - u}^2, \\
S_\mathcal{C}(y,u) &= (L-m) \bigl( u^\tp y - L\norm{y}^2 \bigr) + \norm{u - Ly}^2.
\end{align}
\end{subequations}
This approach of proving that certain quadratic inequalities hold when others do is closely related to the \emph{S-procedure} from control theory~\cite[pp.~23,~33]{lmibook}. The lossless S-procedure states that the two following statements are equivalent:
\begin{enumerate}
	\item For all $x$, if $x^\tp S x \leq 0$, then $x^\tp Q x \leq 0$.\label{Sit1}
	\item There exists $\lambda \geq 0$ such that $S \succeq \lambda Q$.\label{Sit2}
\end{enumerate}
The S-procedure allows inequalities such as~\eqref{eq:sproc_identity_all} to be generated systematically. For example, to generate~\eqref{eq:sproc_identity}, seek a $\lambda \geq 0$ such that
\begin{equation}\label{eq:sp1}
S_\mathcal{C}(y,u) \geq \lambda \left( m^2 \norm{y} - \norm{u}^2 \right) \quad\text{for all }y,u
\end{equation}
The inequality~\eqref{eq:sp1} is equivalent to the semidefinite program:
\begin{equation}\label{eq:sp2}
	\text{Find }\lambda \geq 0\text{ such that: }
	\bmat{ mL & -\frac{L+m}{2} \\ -\frac{L+m}{2} & 1} \succeq \lambda \bmat{ m^2 & 0 \\ 0 & -1}.
\end{equation}
It can be shown that the unique solution to~\eqref{eq:sp2} is $\lambda = \frac{L-m}{2m}$, which leads to~\eqref{eq:sproc_identity}. Since the S-procedure is \emph{lossless} (necessary and sufficient), it will always find an algebraic identity if one exists. In other words, if there is no $\lambda$ that satisfies~\eqref{eq:sp2}, then $S_\mathcal{C}(y,u) \leq 0$ \emph{does not} imply that $m\norm{y} \leq \norm{u}$ for all $y,u$.

\end{document}